\documentclass[12pt]{article}

\usepackage{amsthm}
\usepackage{amssymb}
\usepackage{amsmath}
\usepackage{amsfonts}
\usepackage{enumerate}
\usepackage{verbatim}
\usepackage[dvips]{graphicx}
\usepackage[english]{babel}
\usepackage[cp1250]{inputenc}
\usepackage[T1]{fontenc}
\usepackage[a4paper,left=3cm,right=3cm,top=3cm,bottom=3cm]{geometry}

\newtheorem{col}{Corollary}[section]
\newtheorem{thm}{Theorem}[section]
\newtheorem{lem}{Lemma}[section]

\newtheorem{df}{Definition}[section]

\numberwithin{equation}{section}
\newcommand{\Proof}{\noindent\textbf{{Proof:}}\newline}
\newcommand{\cbdu}{\quad\hfill\mbox{$\Box$}\\[3mm]}

\begin{document}

\title{\textbf{Backward stochastic differential equations with time delayed generators - results and counterexamples}}

\author{\textbf{{\L}ukasz Delong}$^{1}$, \textbf{Peter Imkeller}$^{2}$\\
\
\\
\footnotesize{$^{1}$Institute of Econometrics, Division of Probabilistic Methods}\\
\footnotesize{Warsaw School of Economics}\\
\footnotesize{Al. Niepodleglosci 162, 02-554 Warsaw, Poland}\\
\
\\
\footnotesize{$^{2}$Institut f\"ur Mathematik}\\
\footnotesize{Humboldt-Universit\"at zu Berlin}\\
\footnotesize{Unter den Linden 6, 10099 Berlin, Germany}}

\date{}

\maketitle

\newpage

\begin{abstract}
\noindent We deal with backward stochastic differential equations
with time delayed generators. In this new type of equations, a
generator at time $t$ can depend on the values of a solution in the
past, weighted with a time delay function for instance of the moving
average type. We prove existence and uniqueness of a solution for a
sufficiently small time horizon or for a sufficiently small
Lipschitz constant of a generator. We give examples of BSDE with
time delayed generators that have multiple solutions or that have no
solutions. We show for some special class of generators that
existence and uniqueness may still hold for an arbitrary time
horizon and for arbitrary Lipschitz constant. This class includes
linear time delayed generators, which we study in more detail. We
are concerned with different properties of a solution of a BSDE with
time delayed generator, including the inheritance of boundedness
from the terminal condition, the comparison principle, the existence
of a measure solution and the \emph{BMO} martingale property. We
give examples in which they may fail.

\noindent \textbf{Keywords:} backward stochastic differential
equation, time delayed generator, contraction inequality, comparison
principle, measure solution, \emph{BMO} martingale.

\end{abstract}

\newpage

\section{Introduction}
\indent Backward stochastic differential equations have been
introduced in \cite{PP}. Since then, they have been thoroughly
studied in the literature, see \cite{K} or \cite{I1} and references
therein. The classical theory of BSDE driven by Brownian motions and
with Lipschitz continuous generators has been extended in different
directions. For instance, \cite{Kob} discusses the existence of a
solution in case the generator is of quadratic growth in the control
variable; the existence of a solution for BSDE driven on a more
general stochastic basis, created by L\'{e}vy processes resp.
continuous martingales, is considered respectively in \cite{Bech}
and in \cite{Mor}; a theory of BSDE with random
time horizon is investigated in \cite{BS}. \\
\indent In this paper we study a new class of backward stochastic
differential equation, the dynamics of which is given by
$$Y(t)=\xi+\int_{t}^{T}f(s,Y_{s},Z_{s})ds-\int_{t}^{T}Z(s)dW(s),\quad t\in[0,T].$$
Here, a generator $f$ at time $s$ depends arbitrary on the past
values of a solution $(Y_{s},Z_{s})=(Y(s+u),Z(s+u))_{-T\leq u\leq
0}$. They can be called backward stochastic differential equations
with time delayed generators.  This type of equations has been
investigated for the first time very recently in \cite{BI}, where
only a special form of time delay in $Z$ is considered, namely
$f(s,y_{s},z_{s})=\int_{-T}^{0}g(s+u,z(s+u))\alpha(du)$ with a
measure $\alpha$. The authors prove that in this case, there exists
a unique solution on $[0,T]$ for $T=1$.\\
\indent We aim at providing some contributions to a general theory
of BSDE in which the time delayed generators satisfy Lipschitz
conditions. We are interested in existence and uniqueness results
and in properties of solutions. We prove that a unique solution
exists provided that the generator's Lipschitz constant is
sufficiently small or the evolution is constrained to a sufficiently
small time horizon. The result is optimal as, in general, existence
and uniqueness for a solution of a backward stochastic differential
equation with time delayed (Lipschitz) generator may fail, contrary
to the classical theory of BSDEs without delays, where global
existence and uniqueness results are proved, see \cite{K} or
\cite{I1}. For cases of general Lipschitz constants or time horizons
we give examples of BSDE that have multiple solutions or no
solutions at all. Following \cite{BI}, we also study BSDE with time
delayed generators independent of $y$, and fulfilling Lipschitz
conditions. We show that a unique solution exists if the delay
measure $\alpha$ is supported on $[-\gamma,0]$ with a sufficiently
small time delay $\gamma$. Moreover, in the case of a linear time
delayed generator, which fits into the framework of \cite{BI},
we derive an explicit solution to our BSDE.\\
\indent We further consider properties of solutions of time delayed
BSDE, such as the inheritance of boundedness from the terminal
condition, the comparison principle, measure solutions and the
\emph{BMO} martingale property. All these concepts have turned out
to be very useful in the theory of BSDE without delay, see \cite{I1}
and \cite{K}. We find again that without requiring additional
assumptions these well-known properties, which hold in the classical
setting, may fail for a solution of a time delayed BSDE. We are only
able to prove that the \emph{BMO} martingale property holds in the
case of linear
time delayed generators independent of $Y$. \\
\indent We would like to point out that except \cite{BI} the only
paper we are aware of that deals with BSDE with time delayed
generators is \cite{FT}. In \cite{FT} a forward-backward system of
stochastic differential equations is considered in which the time
delay appears in the forward component, and not in the backward one.
This setting is completely different from the one considered in
\cite{BI} and here. We would like to recall that forward stochastic
differential equations with time delays, called functional
stochastic differential equations, have been studied extensively in
the literature. See for example \cite{Moh}, \cite{Xu}
and references therein.\\
\indent Finally, we would like to refer the reader to the
accompanying paper \cite{DI}, where existence and uniqueness of a
solution of a BSDE driven by a Brownian motion and a Poisson random
measure and with time delayed generator is discussed, together with
its Malliavin's differentiability, both with respect to the
continuous
as well as the jump component.\\
\indent This paper is structured as follows. In Section 2, we deal
with uniqueness and existence of a solution of a backward stochastic
differential equation with time delayed generator. Counterexamples
showing that we cannot obtain unique solutions in a more general
setting are given in Section 3. Linear time delayed generators,
depending only on the control variable $z$, are studied in Section
4, together with the inheritance of boundedness from the terminal
condition and the \emph{BMO} property. Section 5 investigates the
concepts of measure solution and the comparison principle.

\section{Existence and uniqueness of a solution}
\indent We consider a probability space
$(\Omega,\mathcal{F},\mathbb{P})$ with a filtration
$\mathbb{F}=(\mathcal{F}_{t})_{0\leq t\leq T}$ and a finite time
horizon $T<\infty$. We assume that the filtration $\mathbb{F}$ is
the natural filtration generated by a Brownian motion
$W:=(W(t),0\leq t\leq T)$, augmented by all sets of
$\mathbb{P}$-measure zero. \\
\indent We shall work with the following topological vector spaces.
\begin{df}
\begin{itemize}
    \item[1.] Let $L^{2}_{-T}(\mathbb{R})$ denote the space of measurable
        functions $z:[-T,0]\rightarrow\mathbb{R}$ satisfying
        $$\int_{-T}^{0}|z(t)|^{2}dt<\infty.$$
    \item[2.] Let $L^{\infty}_{-T}(\mathbb{R})$ denote the space of bounded, measurable functions
        $y:[-T,0]\rightarrow\mathbb{R}$ satisfying
        $$\sup_{t\in[-T,0]} \big|y(t)\big|^{2}<\infty.$$
    \item[3.] For $p\geq 2$, let $\mathbb{L}^{p}(\mathbb{R})$ denote
        the space of $\mathcal{F}_{T}$-measurable random variables $\xi:\Omega\rightarrow\mathbb{R}$ satisfying
        $$\mathbb{E}\big[\big|\xi\big|^{p}\big]<\infty.$$
    \item [4.]Let $\mathbb{H}^{2}_{T}(\mathbb{R)}$ denote
        the space of $\mathbb{F}$-predictable processes $Z:\Omega\times[0,T]\rightarrow\mathbb{R}$ satisfying
        $$\mathbb{E}\big[\int_{0}^{T}\big|Z(t)\big|^{2}dt\big]<\infty.$$
    \item [5.] Finally, let $\mathbb{S}^{2}_{T}(\mathbb{R)}$ denote
        the space of $\mathbb{F}$-adapted, product measurable processes
        $Y:\Omega\times[0,T]\rightarrow\mathbb{R}$ satisfying
        $$\mathbb{E}\big[\sup_{t\in[0,T]} \big|Y(t)\big|^{2}\big]<\infty.$$
  \end{itemize}
The spaces $\mathbb{H}^{2}_{T}(\mathbb{R)}$ and
$\mathbb{S}^{2}_{T}(\mathbb{R)}$ are endowed with the norms
$$\big\|Z\big\|^{2}_{\mathbb{H}^{2}_{T}}=\mathbb{E}\big[\int_{0}^{T}e^{\beta t}\big|Z(t)\big|^{2}dt\big],$$
$$\big\|Y\big\|^{2}_{\mathbb{S}^{2}_{T}}=\mathbb{E}\big[\sup_{t\in[0,T]} e^{\beta t}\big|Y(t)\big|^{2}\big],$$
with some $\beta>0$.
\end{df}
As usual, by $\lambda$ we denote Lebesgue measure on
$([-T,0],\mathcal{B}([-T,0]))$, where $\mathcal{B}([-T,0])$ stands
for the Borel sets of $[-T,0]$. In the sequel let us simply write
$\mathbb{S}^{2}(\mathbb{R})\times\mathbb{H}^{2}(\mathbb{R})$ for
$\mathbb{S}^{2}_{T}(\mathbb{R})\times\mathbb{H}^{2}_{T}(\mathbb{R})$.\\
\indent We shall deal with the existence and uniqueness of a
solution $(Y,Z):=(Y(t),Z(t))_{0\leq t\leq T}$ of a backward
stochastic differential equation with time delayed generator, the
dynamics of which is given by
\begin{eqnarray}\label{bsde}
Y(t)=\xi+\int_{t}^{T}f(s,Y_{s},Z_{s})ds-\int_{t}^{T}Z(s)dW(s),\quad
0\leq t\leq T,
\end{eqnarray}
where the generator $f$ at time $s\in[0,T]$ depends on the past
values of the solution denoted by $Y_{s}:=(Y(s+u))_{-T\leq u\leq 0}$
and $Z_{s}:=(Z(s+u))_{-T\leq u\leq 0}$. We always set $Z(t)=0$ and $Y(t)=Y(0)$ for $t<0$.\\
\indent We investigate {\eqref{bsde} under the following
assumptions:
\begin{description}
  \item[(A1)] $\xi\in\mathbb{L}^{2}(\mathbb{R})$ for the terminal variable $\xi$,
  \item[(A2)] the generator
  $f:\Omega\times [0,T]\times L^{\infty}_{-T}(\mathbb{R}) \times L^{2}_{-T}(\mathbb{R})\rightarrow \mathbb{R}$
  is product measurable, $\mathbb{F}$-adapted and Lipschitz
  continuous in the sense that for some probability measure $\alpha$ on
  $([-T,0]\times\mathcal{B}([-T,0]))$
  \begin{eqnarray*}
  \lefteqn{|f(t,y_{t},z_{t})-f(t,\tilde{y}_{t},\tilde{z}_{t})|^{2}}\nonumber\\
  &\leq& K\large(\int_{-T}^{0}|y(t+u)-\tilde{y}(t+u)|^{2}\alpha(du)\\
  &&+\int_{-T}^{0}|z(t+u)-\tilde{z}(t+u)|^{2}\alpha(du)\large),
  \end{eqnarray*}
  holds for $\mathbb{P}\times\lambda$-a.e. $(\omega,t)\in\Omega\times[0,T]$ and for any
  $(y_{t},z_{t}),(\tilde{y}_{t},\tilde{z}_{t})\in L^{\infty}_{-T}(\mathbb{R})\times L^{2}_{-T}(\mathbb{R})$.
  \item[(A3)] $\mathbb{E}\big[\int_{0}^{T}|f(t,0,0)|^{2}dt\big]<\infty$,
  \item[(A4)] $f(t,.,.)=0$ for $t<0$.
\end{description}
We remark that $f(t,0,0)$ in \textbf{(A3)} should be understood as a
value of the generator $f(t,y_{t},z_{t})$ at $y(t+u)=z(t+u)=0,-T\leq
u\leq 0$. We would like to point out that the assumption
\textbf{(A4)} in fact allows us to take $Y(t)=Y(0)$ and $Z(t)=0$ for
$t<0$, as a solution of \eqref{bsde}. Examples of generators could
be linear functions of the form
$f(t,y_{t},z_{t})=K\int_{0}^{t}z(s)ds$ or $f(t,y_{t},z_{t})=Kz(t-r),
0\le t\le T$ with a fixed time delay $r$, as studied in more detail
in
Section 4.\\
\indent Note that for
$(Y,Z)\in\mathbb{S}^{2}(\mathbb{R})\times\mathbb{H}^{2}(\mathbb{R})$
the generator is well-defined and $\mathbb{P}$-a.s integrable as
\begin{eqnarray}\label{estf}
\lefteqn{\int_{0}^{T}|f(t,Y_{t},Z_{t})|^{2}dt\leq2\int_{0}^{T}|f(t,0,0)|^{2}dt}\nonumber\\
&&+2K\large(\int_{0}^{T}\int_{-T}^{0}|Y(t+u)|^{2}\alpha(du)dt+\int_{0}^{T}\int_{-T}^{0}|Z(t+u)|^{2}\alpha(du)dt\large)\nonumber\\
&=&2\int_{0}^{T}|f(t,0,0)|^{2}dt\nonumber\\
&&+2K\int_{-T}^{0}\int_{u}^{T+u}|Y(v)|^{2}dv\alpha(du)+2K\int_{-T}^{0}\int_{u}^{T+u}|Z(v)|^{2}dv\alpha(du)\nonumber\\
&\leq
&2\int_{0}^{T}|f(t,0,0)|^{2}dt+2K\big(T\sup_{v\in[0,T]}|Y(v)|^{2}+\int_{0}^{T}|Z(v)|^{2}dv\big)<\infty.
\end{eqnarray}
To justify this, we apply Fubini's theorem, change the variables,
use the assumption that $Z(t)=0$
and $Y(t)=Y(0)$ for $t<0$ and the fact that the probability measure $\alpha$ integrates to $1$.\\
\indent We first state some a priori estimates.
\begin{lem}\label{priori}
Let
$(Y,Z)\in\mathbb{S}^{2}(\mathbb{R})\times\mathbb{H}^{2}(\mathbb{R})$
and
$(\tilde{Y},\tilde{Z})\in\mathbb{S}^{2}(\mathbb{R})\times\mathbb{H}^{2}(\mathbb{R})$
denote solutions of \eqref{bsde} with corresponding parameters
$(\xi,f)$ and $(\tilde{\xi},\tilde{f})$ which satisfy the
assumptions \textbf{(A1)}-\textbf{(A4)}. Then the following
inequalities hold
\begin{eqnarray}\label{p1}
\lefteqn{\|Z-\tilde{Z}\|^{2}_{\mathbb{H}^{2}}}\nonumber\\ &\leq&
e^{\beta T}\mathbb{E}\big[\big|\xi-\tilde{\xi}\big|^{2}\big]+
\frac{1}{\beta}\mathbb{E}\big[\int_{0}^{T}e^{\beta
t}|f(t,Y_{t},Z_{t})-\tilde{f}(t,\tilde{Y}_{t},\tilde{Z}_{t})|^{2}dt\big],
\end{eqnarray}
\begin{eqnarray}\label{p2}
\lefteqn{\|Y-\tilde{Y}\|^{2}_{\mathbb{S}^{2}}}\nonumber\\
&\leq& 8e^{\beta
T}\mathbb{E}\big[\big|\xi-\tilde{\xi}\big|^{2}\big]+
8T\mathbb{E}\big[\int_{0}^{T}e^{\beta
t}|f(t,Y_{t},Z_{t})-\tilde{f}(t,\tilde{Y}_{t},\tilde{Z}_{t})|^{2}dt\big].
\end{eqnarray}
\end{lem}
\Proof The a priori estimates are classical in the theory of BSDEs
without time delay and can be extended to our setting. The
inequality \eqref{p1} follows from Lemma 3.2.1 in \cite{I1}. By
applying It\^{o}'s formula to $e^{\beta t}|Y(t)-\tilde{Y}(t)|^{2}$
on $[0,T]$, taking the expected value and reordering the terms we
derive
\begin{eqnarray*}
\lefteqn{|Y(0)-\tilde{Y}(0)|^{2}+\beta\mathbb{E}\big[\int_{0}^{T}e^{\beta
t}|Y(t)-\tilde{Y}(t)|^{2}dt\big]+
\mathbb{E}\big[\int_{0}^{T}e^{\beta t}|Z(t)-\tilde{Z}(t)|^{2}dt\big]}\\
&\leq&\mathbb{E}\big[e^{\beta
T}|\xi-\tilde{\xi}|^{2}\big]+2\mathbb{E}\big[\int_{0}^{T}e^{\beta
t}|Y(t)-\tilde{Y}(t)||f(t,Y_t,
Z_t)-\tilde{f}(s,\tilde{Y}_t,\tilde{Z}_t)|dt\big].
\end{eqnarray*}
By noticing that
\begin{eqnarray*}
\lefteqn{2|Y(t)-\tilde{Y}(t)||f(t,Y_t,
Z_t)-\tilde{f}(s,\tilde{Y}_t,\tilde{Z}_t)|}\\
&&\leq \beta|Y(t)-\tilde{Y}(t)|^{2}+\frac{1}{\beta}|f(t,Y_t,
Z_t)-\tilde{f}(s,\tilde{Y}_t,\tilde{Z}_t)|^{2},\quad
\mathbb{P}-a.s.,
\end{eqnarray*}
we obtain \eqref{p1}. In order to prove the second inequality, first
notice that for all $t\in[0,T]$
\begin{eqnarray*}
\lefteqn{e^{\frac{\beta}{2}t}|Y(t)-\tilde{Y}(t)|}\\
&\leq& e^{\frac{\beta}{2}T}\mathbb{E}\big[\big|\xi-\tilde{\xi}\big|
|\mathcal{F}_{t}\big]+ \mathbb{E}\big[\int_{0}^{T}e^{\frac{\beta}{2}
s}|f(s,Y_{s},Z_{s})-\tilde{f}(s,\tilde{Y}_{s},\tilde{Z}_{s})|ds|\mathcal{F}_{t}\big].
\end{eqnarray*}
Applying Doob's martingale inequality and Cauchy-Schwarz' inequality
provides the second estimate. \cbdu \\
\indent We state the main theorem of this section.

\begin{thm}\label{thmain}
Assume that \textbf{(A1)}-\textbf{(A4)} hold. For a sufficiently
small time horizon $T$ or for a sufficiently small Lipschitz
constant $K$, the backward stochastic differential equation
\eqref{bsde} has a unique solution
$(Y,Z)\in\mathbb{S}^{2}(\mathbb{R})\times
\mathbb{H}^{2}(\mathbb{R})$.
\end{thm}
\Proof To prove existence and uniqueness of a solution, we follow
the classical idea by constructing a Picard scheme and show its
convergence. See Theorem 2.1
in \cite{K} or Theorem 3.2.1 in \cite{I1}.\\
\indent Let $Y^{0}(t)=Z^{0}(t)=0$ and define recursively for
$n\in\mathbb{N}$
\begin{eqnarray}\label{bsden}
\lefteqn{Y^{n+1}(t)=\xi}\nonumber\\
&&+\int_{t}^{T}f(s,Y^{n}_{s},Z^{n}_{s})ds-\int_{t}^{T}Z^{n+1}(s)dW(s)\quad
0\leq t\leq T.
\end{eqnarray}
\noindent Step 1) Given
$(Y^{n},Z^{n})\in\mathbb{S}^{2}(\mathbb{R})\times
\mathbb{H}^{2}(\mathbb{R})$, the equation \eqref{bsden} has a unique
solution $(Y^{n+1},Z^{n+1})\in\mathbb{S}^{2}(\mathbb{R})\times
\mathbb{H}^{2}(\mathbb{R})$.\\
\noindent Based on the inequality \eqref{estf}, we can conclude that
\begin{eqnarray*}
\lefteqn{\mathbb{E}\big[\int_{0}^{T}|f(t,Y^{n}_{t},Z^{n}_{t})|^{2}dt\big]}\\
&\leq&\mathbb{E}\big[\int_{0}^{T}|f(t,0,0)|^{2}dt\big]
+2K\big(T\|Y^{n}\|_{\mathbb{S}^{2}}+\|Z^{n}\|_{\mathbb{H}^{2}}\big)<\infty.
\end{eqnarray*}
As in the case of a BSDE without a time delay, the martingale
representation provides a unique process
$Z^{n+1}\in\mathbb{H}^{2}(\mathbb{R})$ such that
\begin{eqnarray*}
\xi+\int_{0}^{T}f(t,Y^{n}_{t},Z^{n}_{t})dt=\mathbb{E}\big[\xi+\int_{0}^{T}f(t,Y^{n}_{t},Z^{n}_{t})dt\big]
+\int_{0}^{T}Z^{n+1}(t) dW(t),
\end{eqnarray*}
and we take $Y^{n+1}$ as a continuous version of
\begin{eqnarray*}
Y^{n+1}(t) =
\mathbb{E}\big[\xi|\mathcal{F}_{t}\big]+\mathbb{E}\big[\int_{t}^{T}f(s,Y^{n}_{s},Z^{n}_{s})ds|\mathcal{F}_{t}\big],\quad
0\leq t\leq T.
\end{eqnarray*}
Similarly, as in Lemma \ref{priori}, Doob's inequality,
Cauchy-Schwarz' inequality and the estimates \eqref{estf} yield that
$Y^{n+1}\in\mathbb{S}^{2}(\mathbb{R})$.\\
\noindent In Step 2) we show the convergence of the sequence
$(Y^{n},Z^{n})$ in
$\mathbb{S}^{2}(\mathbb{R})\times\mathbb{H}^{2}(\mathbb{R})$.\\
\noindent The estimates \eqref{p1} and \eqref{p2} give the
inequality
\begin{eqnarray}\label{contr1}
\lefteqn{\big\|Y^{n+1}-Y^{n}\big\|^{2}_{\mathbb{S}^{2}}+\big\|Z^{n+1}-Z^{n}\big\|^{2}_{\mathbb{H}^{2}}}\nonumber\\
&&\leq(8T+\frac{1}{\beta})\mathbb{E}\big[\int_{0}^{T}e^{\beta
t}|f(t,Y^{n}_{t},Z^{n}_{t})-f(t,Y^{n-1}_{t},Z^{n-1}_{t})|^{2}dt\big].
\end{eqnarray}
By applying the Lipschitz condition \textbf{(A2)}, Fubini's theorem,
changing the variables and using the assumption that
$Y^{n}(s)=Y^{n}(0)$ and $Z^{n}(s)=0$ for $s<0$ and all $n\geq 0$, we
can derive
\begin{eqnarray}\label{contr2}
\lefteqn{\mathbb{E}\big[\int_{0}^{T}e^{\beta t}|f(t,Y^{n}_{t},Z^{n}_{t})-f(t,Y^{n-1}_{t},Z^{n-1}_{t})|^{2}dt\big]}\nonumber\\
&\leq& K\mathbb{E}\big[\int_{0}^{T}e^{\beta
t}\int_{-T}^{0}|Y^n(t+u)-Y^{n-1}(t+u)|^{2}\alpha(du)dt\nonumber\\
&&+\int_{0}^{T}e^{\beta t}\int_{-T}^{0}|Z^n(t+u)-Z^{n-1}(t+u)|^{2}\alpha(du)dt\big]\nonumber\\
&=&K\mathbb{E}\big[\int_{-T}^{0}e^{-\beta u}\int_{0}^{T}e^{\beta (t+u)}|Y^n(t+u)-Y^{n-1}(t+u)|^{2}dt\alpha(du)\nonumber\\
&&+\int_{-T}^{0}e^{-\beta u}\int_{0}^{T}e^{\beta (t+u)}|Z^n(t+u)-Z^{n-1}(t+u)|^{2}dt\alpha(du)\big]\nonumber\\
&=&K\mathbb{E}\big[\int_{-T}^{0}e^{-\beta u}\int_{u}^{T+u}e^{\beta v}|Y^n(v)-Y^{n-1}(v)|^{2}dv\alpha(du)\nonumber\\
&&+\int_{-T}^{0}e^{-\beta u}\int_{u}^{T+u}e^{\beta v}|Z^n(v)-Z^{n-1}(v)|^{2}dv\alpha(du)\big]\nonumber\\
&\leq&K\int_{-T}^{0}e^{-\beta
u}\alpha(du)\big(T\big\|Y^{n}-Y^{n-1}\big\|^{2}_{\mathbb{S}^{2}}+\big\|Z^{n}-Z^{n-1}\big\|^{2}_{\mathbb{H}^{2}}\big).
\end{eqnarray}
From \eqref{contr1} and \eqref{contr2}, we obtain
\begin{eqnarray}\label{contr}
\lefteqn{\big\|Y^{n+1}-Y^{n}\big\|^{2}_{\mathbb{S}^{2}}+\big\|Z^{n+1}-Z^{n}\big\|^{2}_{\mathbb{H}^{2}}}\nonumber\\
&&\leq
\delta(T,K,\beta,\alpha)\big(\big\|Y^{n}-Y^{n-1}\big\|^{2}_{\mathbb{S}^{2}}+\big\|Z^{n}-Z^{n-1}\big\|^{2}_{\mathbb{H}^{2}}),
\end{eqnarray}
with
$$\delta(T,K,\beta,\alpha)=(8T+\frac{1}{\beta})K\int_{-T}^{0}e^{-\beta u}\alpha(du)\max\{1,T\}.$$
For $\beta=\frac{1}{T}$ we have
\begin{eqnarray*}
\delta(T,K,\beta,\alpha)\leq 9TKe\max\{1,T\}.
\end{eqnarray*}
For a sufficiently small $T$ or for a sufficiently small $K$, the
inequality \eqref{contr} is a contraction, and there exists a unique
limit
$(Y,Z)\in\mathbb{S}^{2}(\mathbb{R})\times\mathbb{H}^{2}(\mathbb{R})$
of a converging sequence $(Y^{n},Z^{n})_{n\in\mathbb{N}}$, which
satisfies the fixed point equation
\begin{eqnarray*}
Y(t)=
\mathbb{E}\big[\xi|\mathcal{F}_{t}\big]+\mathbb{E}\big[\int_{t}^{T}f(s,Y_{s},Z_{s})ds|\mathcal{F}_{t}\big],\quad
0\leq t\leq T.
\end{eqnarray*}
\noindent Step 4) Define a solution $\bar{Y}$ of \eqref{bsde} as a
continuous version of
\begin{eqnarray*}
\bar{Y}(t)=
\mathbb{E}\big[\xi|\mathcal{F}_{t}\big]+\mathbb{E}\big[\int_{t}^{T}f(s,Y_{s},Z_{s})ds|\mathcal{F}_{t}\big],\quad
0\leq t\leq T,
\end{eqnarray*}
where $(Y,Z)$ is the limit constructed in Step 3. \cbdu

\indent Theorem \ref{thmain} triggers the immediate question: is it
possible to obtain existence and uniqueness for a bigger time
horizon $T$ and/or an arbitrary Lipschitz constant $K$? In the
following section we show that such an extension is not possible.
However, for a special class of generators Theorem \ref{thmain} may
be generalized, as we now show.
\begin{thm}
Assume that \textbf{(A1)}-\textbf{(A4)} hold and that the generator
is independent of $y_{t}$, i. e. for $t\in[0,T]$ we have
$f(t,y_{t},z_{t})=f(t,z_{t})$. Let the measure $\alpha$ be supported
on the interval $[-\gamma,0]$, where $\gamma$ is a constant. For a
sufficiently small time delay $\gamma$ the backward stochastic
differential equation \eqref{bsde} has a unique solution
$(Y,Z)\in\mathbb{S}^{2}(\mathbb{R})\times
\mathbb{H}^{2}(\mathbb{R})$.
\end{thm}
\Proof The proof is very similar to the previous one. Note that in
this case, based on \eqref{contr1} and \eqref{contr2}, we have
\begin{eqnarray*}
\big\|Z^{n+1}-Z^{n}\big\|^{2}_{\mathbb{H}^{2}}\leq
\delta(T,K,\beta,\alpha)\big\|Z^{n}-Z^{n-1}\big\|^{2}_{\mathbb{H}^{2}},
\end{eqnarray*}
with
$$\delta(T,K,\beta,\alpha)=\frac{1}{\beta}K\int_{-\gamma}^{0}e^{-\beta u}\alpha(du)\leq
\frac{Ke^{\beta\gamma}}{\beta}\int_{-\gamma}^{0}\alpha(du)=\frac{Ke^{\beta\gamma}}{\beta},$$
which is smaller than $1$ for sufficiently big $\beta$ and small
$\gamma$. This proves the convergence of $(Z^{n})_{n\in\mathbb{N}}$.
To get the convergence of $(Y^{n})_{n\in\mathbb{N}}$, notice again
that, by \eqref{contr1} and \eqref{contr2}
\begin{eqnarray*}
\big\|Y^{n+1}-Y^{n}\big\|^{2}_{\mathbb{S}^{2}}\leq
8T\delta(T,K,\beta,\alpha)\big\|Z^{n}-Z^{n-1}\big\|^{2}_{\mathbb{H}^{2}}.
\end{eqnarray*}
\cbdu

\indent Finally, to complete our presentation of the current state
of knowledge on BSDE with time delayed generators, we shall recall a
theorem proved in \cite{BI} recently.
\begin{thm}\label{thimkeller}
Assume that $\xi\in\mathbb{L}^{2+\epsilon}(\mathbb{R})$ for some
$\epsilon>0$ and that \textbf{(A2)}-\textbf{(A4)} hold with respect
to a generator of the form
\begin{eqnarray*}
f(t,z_{t})=\int_{-T}^{0}g(t+u,z(t+u))\alpha(du),\quad z\in
L^{2}_{-T}(\mathbb{R}),
\end{eqnarray*}
where $\alpha$ is a finite measure. The BSDE \eqref{bsde} has a
unique solution $(Y,Z)\in\mathbb{S}^{2}(\mathbb{R})\times
\mathbb{H}^{2}(\mathbb{R})$.
\end{thm}
We remark that our Theorem \ref{thimkeller} is a slight extension of
the theorem stated in \cite{BI}.  First, it is straightforward to
extend the result of \cite{BI} from $T=1$ to an arbitrary $T$.
Secondly, it is possible to prove the result of \cite{BI} under
weaker integrability assumptions concerning $\xi$, by following the
proof in \cite{BI} and replacing the Cauchy-Schwarz inequality with
the H\"older inequality. This allows for
$\xi\in\mathbb{L}^{2+\epsilon}(\mathbb{R})$ in Theorem
\ref{thimkeller}.

\section{Non-uniqueness and multiple solutions}
\indent In this section we discuss examples of BSDE with time
delayed generators that fail to have solutions or have more than
one. This confirms that there is a natural boundary for extensions
of the local existence and uniqueness result from Theorem
\ref{thmain}, and that one cannot expect to have existence and
uniqueness for an arbitrary time horizon $T$ and an arbitrary
Lipschitz constant $K$
without additional requirements.\\
\indent {\bf Example 1}\\
Let us first investigate the backward stochastic differential
equation with the following generator of Lipschitz constant $K>0$
and of fixed time delay
\begin{eqnarray}\label{exm1}
Y(t)=\xi+\int_{t}^{T}KY(s-T)ds-\int_{t}^{T}Z(s)dW(s),\quad 0\leq
t\leq T.
\end{eqnarray}
Using the notions of the previous section, (\ref{exm1}) can be
rewritten as
\begin{eqnarray*}
\lefteqn{Y(t)=\xi}\\
&&+\int_{t}^{T}\int_{-T}^{0}KY(s+u)\mathbf{1}_{[0,\infty)}(s)\,\alpha(du)ds-\int_{t}^{T}Z(s)dW(s),\quad
0 \leq t\leq T,
\end{eqnarray*}
with Dirac measure $\alpha$ concentrated at the point $T$. The
equation \eqref{exm1} is clearly equivalent to
\begin{eqnarray}\label{exm2}
Y(t)=\xi+K(T-t)Y(0)-\int_{t}^{T}Z(s)dW(s),\quad 0\leq t\leq T.
\end{eqnarray}
For $t=0$ we arrive at
\begin{eqnarray*}
(1-TK)Y(0)=\xi-\int_{0}^{T}Z(s)dW(s),
\end{eqnarray*}
and integrating on both sides produces the condition
\begin{eqnarray*}\label{cond}
\mathbb{E}[\xi]=(1-TK)Y(0).
\end{eqnarray*}
We consider three cases.\\
\noindent Case 1) $TK< 1$.\\
\noindent Define $Z$ as the unique square integrable process from
the martingale representation of $\xi\in\mathbb{L}^{2}(\mathbb{R})$
given by
\begin{eqnarray*}
\xi=\mathbb{E}[\xi]+\int_{0}^{T}Z(s)dW(s),
\end{eqnarray*}
and the process $Y$, according to \eqref{exm2}, by
\begin{eqnarray}
Y(t)&=&\mathbb{E}[\xi]+K(T-t)Y(0)+\int_{0}^{t}Z(s)dW(s)\nonumber\\
&=&Y(0)(1-tK)+ \int_{0}^{t}Z(s)dW(s)\nonumber\\
&=&\frac{1-tK}{1-TK}\mathbb{E}[\xi]+ \int_{0}^{t}Z(s)dW(s),\quad
0\leq t\leq T.
\end{eqnarray}
The pair $(Y,Z)\in\mathbb{S}^{2}(\mathbb{R})\times
\mathbb{H}^{2}(\mathbb{R})$ is the unique solution of \eqref{exm1}
on $[0,T]$. In fact, suppose there were another solution
$(\tilde{Y},\tilde{Z})$ of \eqref{exm1} on $[0,T]$. Writing the
difference of \eqref{exm2} for the two solutions, and using
$(1-TK)Y(0)=\mathbb{E}[\xi]=(1-TK)\tilde{Y}(0)$ for getting $Y(0) =
\tilde{Y}(0)$, we obtain
$$\int_0^T (Z(s) - \tilde{Z}(s)) d W(s) =0,\quad \mathbb{P}-a.s.,$$
hence $Z = \tilde{Z}$, whence finally $Y = \tilde{Y}.$ \\
\noindent Case 2) $TK=1$ and $\mathbb{E}[\xi]\neq 0$. \\
\noindent The condition $\mathbb{E}[\xi]=(1-TK)Y(0)$ is not
satisfied and therefore equation \eqref{exm1}
does not have any solution.\\
\noindent Case 3) $TK=1$ and $\mathbb{E}[\xi]=0$.\\
As in case 1), define $Z$ as the unique square integrable process
appearing in the martingale representation of
$\xi\in\mathbb{L}^{2}(\mathbb{R})$, and the process $Y$ as
\begin{eqnarray}\label{exm3.1}
Y(t)=Y(0)(1-tK)+ \int_{0}^{t}Z(s)dW(s),\quad 0\leq t\leq T,
\end{eqnarray}
with an arbitrary $Y(0)\in\mathbb{L}^{2}(\mathbb{R})$ which is
$\mathcal{F}_{0}$-measurable. Any pair
$(Y,Z)\in\mathbb{S}^{2}(\mathbb{R})\times
\mathbb{H}^{2}(\mathbb{R})$ satisfying \eqref{exm3.1} is a solution
of \eqref{exm1} on $[0,T]$.\\\par\bigskip\bigskip

\indent {\bf Example 2}\\
Next, again let $K\in\mathbb{R}$, so that $|K|$ stands for the
Lipschitz constant of the time delayed generator, we study the
backward stochastic differential equation
\begin{eqnarray}\label{exm4}
Y(t)=\xi+\int_{t}^{T}\int_{0}^{s}KY(u)duds-\int_{t}^{T}Z(s)dW(s),\quad
0\leq t\leq T.
\end{eqnarray}
With the notation of the previous section the equation is of the
form
\begin{eqnarray*}
\lefteqn{Y(t)=\xi}\\
&&+\int_{t}^{T}\int_{-T}^{0}KTY(s+u)\mathbf{1}\{s+u\geq
0\}\alpha(du)ds\\
&&-\int_{t}^{T}Z(s)dW(s),\quad 0\leq t\leq T,
\end{eqnarray*}
with a uniform measure $\alpha$ on $[-T,0]$. Changing the order of
integration in the generator term and calculating the difference
$Y(t)-Y(0)$ yields
\begin{eqnarray}\label{exm44}
\lefteqn{Y(t)=Y(0)}\nonumber\\
&&-K\int_{0}^{t}(t-s)Y(s)ds+\int_{0}^{t}Z(s)dW(s),\quad 0\leq t\leq
T.
\end{eqnarray}
In the sequel we construct a solution of \eqref{exm4}. We comment on
the main steps and leave details of the tedious but simple
calculations to
the reader.\\
\indent Consider for a moment the deterministic integral equation
corresponding to \eqref{exm44}
\begin{eqnarray}\label{exm4int}
y(t)=y(0)-K\int_{0}^{t}(t-s)y(s)ds+h(t),
\end{eqnarray}
with a twice continuously differentiable function
$h\in\mathcal{C}^{2}(\mathbb{R})$ such that $h(0)=0$ and with a
given initial condition $y(0)$. By differentiating, we obtain the
nonhomogeneous linear second order differential equation
\begin{eqnarray}\label{exm4dif}
y^{''}(t)+Ky(t)=h^{''}(t).
\end{eqnarray}
The fundamental solution of the homogeneous part of \eqref{exm4dif}
is well known and its general form is
\begin{eqnarray*}
y(t)=Ae^{\sqrt{-K}t}+Be^{-\sqrt{-K}t},
\end{eqnarray*}
where $A, B$ are constants, and $\sqrt{-K}$ for $K>0$ is understood
as a complex number. It is easy to check that the following formula
gives a general solution of the inhomogeneous equation
\eqref{exm4dif}:
\begin{eqnarray*}
\lefteqn{y(t)=Ae^{\sqrt{-K}t}+Be^{-\sqrt{-K}t}}\nonumber\\
&&+\int_{0}^{t}\frac{h^{''}(s)e^{-\sqrt{-K}s}}{2\sqrt{-K}}dse^{\sqrt{-K}t}-
\int_{0}^{t}\frac{h^{''}(s)e^{\sqrt{-K}s}}{2\sqrt{-K}}dse^{-\sqrt{-K}t}.
\end{eqnarray*}
Integrating by parts twice gives
\begin{eqnarray}\label{exm4difsol}
\lefteqn{y(t)=Ae^{\sqrt{-K}t}+Be^{-\sqrt{-K}t}+h(t)}\nonumber\\
&&+\int_{0}^{t}h(s)e^{-\sqrt{-K}s}ds\frac{\sqrt{-K}}{2}e^{\sqrt{-K}t}-
\int_{0}^{t}h(s)e^{\sqrt{-K}s}ds\frac{\sqrt{-K}}{2}e^{-\sqrt{-K}t}.
\end{eqnarray}
One can further check that the part of the solution
\eqref{exm4difsol} containing $h$ satisfies the integral equation
\eqref{exm4int} even without any  differentiability assumptions
concerning $h$. One can finally derive the following conditions,
under which \eqref{exm4difsol} solves the integral equation
\eqref{exm4int}:
\begin{eqnarray*}
\left\{\begin{array}{ll}
A=B,\\
A+B=y(0).
\end{array}\right.
\end{eqnarray*}
Returning to our backward stochastic differential equation, it is
straightforward to replace $h$ with $\int_0^\cdot Z(s) d W(s)$, and
to conclude that a solution of \eqref{exm44} must be of the form
\begin{eqnarray*}
\lefteqn{Y(t)=\frac{Y(0)}{2}\big(e^{\sqrt{-K}t}+e^{-\sqrt{-K}t}\big)+\int_{0}^{t}Z(s)dW(s)}\nonumber\\
&&+\int_{0}^{t}\int_{0}^{s}Z(u)dW(u)e^{-\sqrt{-K}s}ds\frac{\sqrt{-K}}{2}e^{\sqrt{-K}t}\nonumber\\
&&-
\int_{0}^{t}\int_{0}^{s}Z(u)dW(u)e^{\sqrt{-K}s}ds\frac{\sqrt{-K}}{2}e^{-\sqrt{-K}t},\quad
0\leq t\leq T.
\end{eqnarray*}
By applying Fubini's theorem for stochastic integrals, see Theorem
4.65 in \cite{Pr}, we finally derive
\begin{eqnarray*}
Y(t)=\frac{Y(0)}{2}\big(e^{\sqrt{-K}t}+e^{-\sqrt{-K}t}\big)
+\frac{1}{2}\int_{0}^{t}Z(s)\big(e^{\sqrt{-K}(t-s)}+e^{-\sqrt{-K}(t-s)}\big)dW(s),
\end{eqnarray*}
for $t\in[0,T]$, and $Y(0)$ is determined by
\begin{eqnarray*}
\mathbb{E}[\xi]=\frac{Y(0)}{2}\big(e^{\sqrt{-K}T}+e^{-\sqrt{-K}T}\big).
\end{eqnarray*}

\noindent Case 1) Let us assume that $K<0$. In this case the unique
solution
$(Y,Z)\in\mathbb{S}^{2}(\mathbb{R})\times\mathbb{H}^{2}(\mathbb{R})$
of \eqref{exm4} is given by
\begin{eqnarray*}
Y(t)&=&\mathbb{E}[\xi]\frac{e^{\sqrt{-K}t}+e^{-\sqrt{-K}t}}{e^{\sqrt{-K}T}+e^{-\sqrt{-K}T}}\\
&&+\frac{1}{2}\int_{0}^{t}Z(s)\big(e^{\sqrt{-K}(t-s)}+e^{-\sqrt{-K}(t-s)}\big)dW(s),\nonumber\\
Z(t)&=&\frac{2M(t)}{e^{\sqrt{-K}(T-t)}+e^{-\sqrt{-K}(T-t)}},
\end{eqnarray*}
for $0\leq t\leq T$, where $M:=(M(t))_{0\leq t\leq T}$ is the unique
square integrable process appearing in the martingale representation
of $\xi\in\mathbb{L}^{2}(\mathbb{R})$, namely
\begin{eqnarray*}
\xi=\mathbb{E}[\xi]+\int_{0}^{T}M(t)dW(t).
\end{eqnarray*}

\noindent Case 2) Let us assume now that $K>0$. This case is more
interesting, as it allows for uniqueness, nonexistence and
multiplicity of solutions. By Euler's formula
\begin{eqnarray*}
Y(t)&=&Y(0)\cos(t\sqrt{K})+\int_{0}^{t}\cos((t-s)\sqrt{K})Z(s)dW(s),\quad
0\leq t\leq T,
\end{eqnarray*}
with
\begin{eqnarray*}
\mathbb{E}[\xi]=Y(0)\cos(T\sqrt{K}).
\end{eqnarray*}
\noindent Case 2.1). $T\sqrt{K}<\frac{\pi}{2}$.\\
\noindent The unique solution
$(Y,Z)\in\mathbb{S}^{2}(\mathbb{R})\times\mathbb{H}^{2}(\mathbb{R})$
of \eqref{exm4} is given by
\begin{eqnarray*}
Y(t)&=&\mathbb{E}[\xi]\frac{\cos(t\sqrt{K})}{\cos(T\sqrt{K})}+\int_{0}^{t}\frac{\cos((t-s)\sqrt{K})}{\cos((T-s)\sqrt{K})}M(s)dW(s),\quad
0\leq t\leq T,
\end{eqnarray*}
where $M$ is the unique process arising from the martingale
representation of
$\xi\in\mathbb{L}^{2}(\mathbb{R})$.\\
\noindent Case 2.2). $T\sqrt{K}=\frac{\pi}{2}$ and
$\mathbb{E}[\xi]\neq0$. \\
\noindent Equation \eqref{exm4} does not have any solution, since
condition $\mathbb{E}[\xi]=Y(0)\cos(T\sqrt{K})$ is not satisfied.\\
\noindent Case 2.3). $T\sqrt{K}=\frac{\pi}{2}$ and
$\mathbb{E}[\xi]=0$.\\
\noindent Equation \eqref{exm4} may not have any solution, or may
have multiple solutions. Consider again the representation
\begin{eqnarray*}
\xi&=&\int_{0}^{T}M(s)dW(s),
\end{eqnarray*}
and put
\begin{eqnarray*}
Z(t)=\frac{M(t)}{\cos((T-t)\sqrt{K})}\mathbf{1}\{t>0\},\quad 0\leq
t\leq T.
\end{eqnarray*}
\noindent Case 2.3.a) If
$\mathbb{E}[\int_{0}^{T}|Z(s)|^{2}ds]=+\infty$, then equation
\eqref{exm4} does not have any solution.\\
\noindent Case 3b) If $\mathbb{E}[\int_{0}^{T}|Z(s)|^{2}ds]<\infty$,
then equation \eqref{exm4} has multiple solutions
$(Y,Z)\in\mathbb{S}^{2}(\mathbb{R})\times\mathbb{H}^{2}(\mathbb{R})$
given by
\begin{eqnarray*}
Y(t)&=&Y(0)\cos(t\sqrt{K})+\int_{0}^{t}\cos((t-s)\sqrt{K})Z(s)dW(s),\quad
0\leq t\leq T,
\end{eqnarray*}
with an arbitrary $Y(0)\in\mathbb{L}^{2}(\mathbb{R})$ which is $\mathcal{F}_{0}$-measurable.\\
\noindent To make the example complete, take $K=1$ and notice that
for $\xi=\int_{0}^{\frac{\pi}{2}}\cos(\frac{\pi}{2}-s)dW(s)$ we have
multiple solutions, whereas for $\xi=W(T)$ we don't have any
solution, since
$\mathbb{E}[\int_{0}^{\frac{\pi}{2}}|\frac{1}{\cos(\frac{\pi}{2}-t)}|^{2}ds]=+\infty$.

\section{BSDEs with linear time delayed generators}
In this section we investigate in more details the following
backward stochastic differential equation with a linear time delayed
generator
\begin{eqnarray}\label{bsdelin}
\lefteqn{Y(t)=\xi}\nonumber\\
&&+\int_{t}^{T}\int_{-T}^{0}g(s+u)
Z(s+u)\alpha(du)ds-\int_{t}^{T}Z(s)dW(s),\quad 0\leq t\leq T,
\end{eqnarray}
with
\begin{description}
\item[(A5)] a measurable, uniformly bounded function
$g:[0,T]\rightarrow\mathbb{R}$ and the assumption $g(t)=0$ for
$t<0$.
\end{description}
As for the measure $\alpha$, we are particularly interested in the
two extreme cases in which $\alpha$ is uniform or
a Dirac measure.\\
\indent For the linear equations \eqref{bsdelin} it is possible to
describe solutions explicitly.
\begin{thm}\label{thmlin}
Assume that $\xi\in\mathbb{L}^{2+\epsilon}(\mathbb{R})$ for some
$\epsilon>0$, and \textbf{(A5)} holds. The backward stochastic
differential equation \eqref{bsdelin} has a unique solution
$(Y,Z)\in\mathbb{S}^{2}(\mathbb{R})\times
\mathbb{H}^{2}(\mathbb{R})$, where $Z$ is the process appearing in
the martingale representation
\begin{eqnarray}\label{repr}
\xi=\mathbb{E}^{\mathbb{Q}}[\xi]+\int_{0}^{T}Z(s)dW^{\mathbb{Q}}(s),
\end{eqnarray}
under the equivalent probability measure $\mathbb{Q}$ given by the
density
\begin{eqnarray}\label{measure}
\lefteqn{\frac{d\mathbb{Q}}{d\mathbb{P}}\big|_{\mathcal{F}_{T}}}\nonumber\\
&&=\exp\big(\int_{0}^{T}\alpha((s-T,0])g(s)dW(s)-
\frac{1}{2}\int_{0}^{T}\alpha^{2}((s-T,0])g^{2}(s)ds\big),
\end{eqnarray}
with a $\mathbb{Q}$-Brownian motion
\begin{eqnarray*}
W^{\mathbb{Q}}(t)=W(t)-\int_{0}^{t}\alpha((s-T,0])g(s)ds,\quad 0\leq
t\leq T,
\end{eqnarray*}
and the process $Y$ defined by
\begin{eqnarray}\label{liny}
Y(t)=\mathbb{E}^{\mathbb{Q}}\big[\xi|\mathcal{F}_{t}]+\int_{0}^{t}\alpha((s-T,s-t])Z(s)ds,\quad
0\leq t\leq T.
\end{eqnarray}
\end{thm}
\Proof First we notice that the generator of the equation
\eqref{bsdelin} is Lipschitz continuous in the sense of
\textbf{(A2)} from Section 2, since for $\lambda$-a.e. $t\in[0,T]$
\begin{eqnarray*}
\lefteqn{\big|\int_{-T}^{0}g(t+u)z(t+u)\alpha(du)-\int_{-T}^{0}g(t+u)\tilde{z}(t+u)\alpha(du)\big|^{2}}\nonumber\\
&\leq&
\int_{-T}^{0}g^{2}(t+u)\alpha(du)\int_{-T}^{0}\big|z(t+u)-\tilde{z}(t+u)\big|^{2}\alpha(du)\nonumber\\
&\leq&G\int_{-T}^{0}\big|z(t+u)-\tilde{z}(t+u)\big|^{2}\alpha(du),
\end{eqnarray*}
where we apply Cauchy-Schwarz' inequality, and $G$ denotes the
uniform bound on $g$. By recalling Theorem \ref{thimkeller} we can
conclude that there exists a unique solution
$(Y,Z)\in\mathbb{S}^{2}(\mathbb{R})\times
\mathbb{H}^{2}(\mathbb{R})$ of \eqref{bsdelin}.\\
\indent Notice that by applying Fubini's theorem, changing the
variables and using the assumption that $g$ vanishes for $t<0$, we
can derive an alternative form of the integral of the generator
\begin{eqnarray*}
\int_{t}^{T}\int_{-T}^{0}g(s+u)
Z(s+u)\alpha(du)ds&=&\int_{-T}^{0}\int_{t}^{T}g(s+u)
Z(s+u)ds\alpha(du)\nonumber\\
&=&\int_{-T}^{0}\int_{(t-u)\vee 0}^{T-u}g(v)
Z(v)dv\alpha(du)\nonumber\\
&=&\int_{0}^{T}\int_{v-T}^{(v-t)\wedge 0}g(v) Z(v)\alpha(du)dv,
\end{eqnarray*}
for $0\leq t\leq T$. This allows us to rewrite the BSDE
\eqref{bsdelin} as
\begin{eqnarray*}
Y(t)&=&\xi+\int_{0}^{T}\alpha((s-T,(s-t)\wedge0])g(s)Z(s)ds-\int_{t}^{T}Z(s)dW(s)\nonumber\\
&=&\xi+\int_{0}^{t}\alpha((s-T,s-t])g(s)Z(s)ds\\
&&-\int_{t}^{T}Z(s)\big(dW(s)-\alpha((s-T,0])g(s)ds\big), \quad
0\leq t\leq T.
\end{eqnarray*}
The measure defined in \eqref{measure} is an equivalent probability
measure because $\int_{0}^{T}\alpha^{2}((s-T,0\rangle)g^{2}(s)ds$ is
finite, and hence Novikov's condition is satisfied. We can therefore
deal with the following equation under the measure $\mathbb{Q}$
\begin{eqnarray}\label{bsdelinq}
\lefteqn{Y(t)=\xi}\nonumber\\
&&+\int_{0}^{t}\alpha((s-T,s-t\rangle)g(s)Z(s)ds-\int_{t}^{T}Z(s)dW^{\mathbb{Q}}(s)\big),\quad
0\leq t\leq T,
\end{eqnarray}
where $W^{\mathbb{Q}}$ is a $\mathbb{Q}$-Brownian motion on $[0,T]$.
The assumption that $\xi\in\mathbb{L}^{2+\epsilon}(\mathbb{R})$, for
some $\epsilon>0$, under the measure $\mathbb{P}$, yields that
$\xi\in\mathbb{L}^{2+\frac{\epsilon}{2}}(\mathbb{R})$ under the
measure $\mathbb{Q}$. This is justified by the inequality
\begin{eqnarray*}
\mathbb{E}^{\mathbb{Q}}\big[|\xi|^{2+\frac{\epsilon}{2}}\big]\leq
\big(\mathbb{E}^{\mathbb{P}}\big[\big(\frac{d\mathbb{Q}}{d\mathbb{P}}\big)^{\frac{4+2\epsilon}{\epsilon}}\big]\big)^{\frac{\epsilon}{4+2\epsilon}}
\big(\mathbb{E}^{\mathbb{P}}\big[|\xi|^{2+\epsilon}\big]\big)^{\frac{2+\frac{\epsilon}{2}}{2+\epsilon}}<\infty,
\end{eqnarray*}
which uses that the density $\frac{d\mathbb{Q}}{d\mathbb{P}}$
possesses moments of all orders. Define $(Y,Z)$ according to
\eqref{repr} and \eqref{liny}. Clearly, $(Y,Z)$ is a solution of
\eqref{bsdelinq}. The martingale representation theorem in
$\mathbb{L}^{2+\frac{\epsilon}{2}}(\mathbb{R})$ under $\mathbb{Q}$
provides $Z$ such that
$\mathbb{E}^{\mathbb{Q}}\big[\big(\int_{0}^{T}|Z(s)|^{2}ds\big)^{1+\frac{\epsilon}{4}}\big]<\infty$,
see Theorem 5.1 in \cite{K}. The process $Z$ is also square
integrable under $\mathbb{P}$, as is seen by
\begin{eqnarray*}
\mathbb{E}^{\mathbb{P}}\big[\int_{0}^{T}|Z(s)|^{2}ds\big]\leq
\big(\mathbb{E}^{\mathbb{Q}}\big[\big(\frac{d\mathbb{P}}{d\mathbb{Q}}\big)^{\frac{4+\epsilon}{\epsilon}}\big]\big)^{\frac{\epsilon}{4+\epsilon}}
\big(\mathbb{E}^{\mathbb{Q}}\big[\big(\int_{0}^{T}|Z(s)|^{2}ds\big)^{1+\frac{\epsilon}{4}}
\big]\big)^{\frac{4}{4+\epsilon}}<\infty.
\end{eqnarray*}
It is standard to prove that $Y\in\mathbb{S}^{2}(\mathbb{R})$. We
can conclude now that $(Y,Z)$ defined by \eqref{repr} and
\eqref{liny} is the unique solution of \eqref{bsdelin}. \cbdu
\indent We can state two corollaries.
\begin{col}\label{propint}
Under the assumptions of Theorem \ref{thmlin}, the equation with a
delay distributed according to a uniform measure
\begin{eqnarray*}
Y(t)=\xi+\int_{t}^{T}\int_{0}^{s}
KZ(u)duds-\int_{t}^{T}Z(s)dW(s),\quad 0\leq t\leq T,
\end{eqnarray*}
has a unique solution $(Y,Z)\in\mathbb{S}^{2}(\mathbb{R})\times
\mathbb{H}^{2}(\mathbb{R})$, with $Z$ given by the martingale
representation
\begin{eqnarray*}
\xi=\mathbb{E}^{\mathbb{Q}}[\xi]+\int_{0}^{T}Z(s)dW^{\mathbb{Q}}(s),
\end{eqnarray*}
under the equivalent probability measure
\begin{eqnarray*}
\frac{d\mathbb{Q}}{d\mathbb{P}}\big|_{\mathcal{F}_{T}}=\exp\big(K\int_{0}^{T}(T-s)dW(s)-
K^{2}\frac{1}{2}\int_{0}^{T}(T-s)^{2}ds\big),
\end{eqnarray*}
and the process $Y$ defined by
\begin{eqnarray*}
Y(t)=\mathbb{E}^{\mathbb{Q}}\big[\xi|\mathcal{F}_{t}]+(T-t)K\int_{0}^{t}Z(s)ds,\quad
0\leq t\leq T.
\end{eqnarray*}
\end{col}
The proof is a straightforward application of Theorem \ref{thmlin}.

\begin{col}
Under the assumptions of Theorem \ref{thmlin}, the equation with a
delay distributed according to Dirac measure at the point
$r\in[0,T]$
\begin{eqnarray*}
Y(t)=\xi+\int_{t}^{T}KZ(s-r)ds-\int_{t}^{T}Z(s)dW(s),\quad 0\leq
t\leq T,
\end{eqnarray*}
has a unique solution $(Y,Z)\in\mathbb{S}^{2}(\mathbb{R})\times
\mathbb{H}^{2}(\mathbb{R})$, given by the following statements.
\begin{description}
  \item[1.] On the interval $[0,T-r)$, define $Z$ as the process arising in the martingale representation
  \begin{eqnarray*}
  \mathbb{E}[\xi|\mathcal{F}_{T-r}]=\mathbb{E}^{\mathbb{Q}}\big[\mathbb{E}[\xi|\mathcal{F}_{r}]\big]+\int_{0}^{T-r}Z(s)dW^{\mathbb{Q}}(s),
  \end{eqnarray*}
  under the equivalent probability measure
  \begin{eqnarray*}
  \frac{d\mathbb{Q}}{d\mathbb{P}}\big|_{\mathcal{F}_{T-r}}=\exp\big(KW(T-r)-\frac{1}{2}K^{2}(T-r)\big),
  \end{eqnarray*}
  and
  \begin{eqnarray*}
  Y(t)=\mathbb{E}^{\mathbb{Q}}\big[\mathbb{E}[\xi|\mathcal{F}_{T-r}]|\mathcal{F}_{t}\big]+K\int_{(t-r)\wedge 0}^{t}
  Z(s)ds,\quad 0\leq t< T-r.
  \end{eqnarray*}
  \item[2.] On the interval $[T-r,T]$, define $Z$ as the process arising in the martingale representation
  \begin{eqnarray*}
  \xi=\mathbb{E}[\xi|\mathcal{F}_{T-r}]+\int_{T-r}^{T}Z(s)dW(s),
  \end{eqnarray*}
  and
  \begin{eqnarray*}
  Y(t)=\mathbb{E}[\xi|\mathcal{F}_{t}]+K\int_{(t-r)\wedge 0}^{T-r}Z(s)ds,\quad T-r\leq t\leq T.
  \end{eqnarray*}
\end{description}
\end{col}
\Proof First we notice that in the case of a Dirac measure
concentrated at $r$ we have
$\alpha((s-T,0])=\mathbf{1}_{[0,T-r)}(s)$. We conclude that
$\alpha((s-T,0])=0$ for $s\in[T-r,T]$ and the $\mathbb{Q}$-Brownian
motion is the $\mathbb{P}$-Brownian motion on the interval $[T-r,T]$
(given $\mathcal{F}_{T-r}$). As $\alpha((s-T,0])=1$ for
$s\in[0,T-r)$ we can define the corresponding density process
$\frac{d\mathbb{Q}}{d\mathbb{P}}\big|_{\mathcal{F}_{t}}$. Moreover,
notice that
$\alpha((s-T,s-t])=\mathbf{1}_{[t-r,T-r)}(s)$.\\
\noindent Consider $t\in[0,T-r)$. By taking the conditional
expectation under $\mathbb{P}$, we derive
\begin{eqnarray*}
Y(t)=\mathbb{E}\big[\xi\big|\mathcal{F}_{T-r}\big]+K\int_{0}^{t}\mathbf{1}_{[t-r,T-r)}(s)\,
Z(s)ds-\int_{t}^{T-r}Z(s)dW^{\mathbb{Q}}(s).
\end{eqnarray*}
Defining $(Y,Z)$ as in the the first part of the statement, we get a solution on $[0,T-r)$.\\
\noindent Consider $t\in[T-r,T]$. We now have to deal with the
equation
\begin{eqnarray*}
Y(t)=\xi+K\int_{0}^{t}\mathbf{1}_{[t-r,T-r)}(s)\,Z(s)ds-\int_{t}^{T}Z(s)dW(s).
\end{eqnarray*}
Defining $(Y,Z)$ as in the second part of the statement, we obtain a
solution on $[T-r,T]$. One can verify, similarly as in the proof of
Theorem \ref{thmlin}, that the solution constructed in this way
belongs to the right space:
$(Y,Z)\in\mathbb{S}^{2}(\mathbb{R})\times
\mathbb{H}^{2}(\mathbb{R})$. This completes the proof.\cbdu

\indent We next consider the \emph{BMO} (bounded mean oscillation)
property of the stochastic integral process $\int ZdW$ of the
control component $Z$ of the solution, which is important for
solutions of BSDE without time delay in the generator. We discuss
the question whether this property continues to hold if the
generator possesses some linear time delay feature, such as in
\eqref{bsdelin}.\\
\indent It is well-known, see Lemma 3.1.2 in \cite{I1}, that for a
BSDE with generator not subject to a time delay, a terminal
condition $\xi$ that is $\mathbb{P}$-a.s. bounded, and satisfying
appropriate assumptions, the integral process $\int ZdW$ is a
\emph{BMO} martingale. It is a rather easy exercise to prove this
result for a linear Lipschitz generator without a time delay. We
show that for linear BSDEs with time delayed generators, this
property still holds.\\
\indent Recall that a stochastic integral process
$(\int_{0}^{t}Z(s)dW(s))_{0\leq t\leq T}$ is a \emph{BMO} martingale
iff
\begin{eqnarray*}
\sup_{\tau}\mathbb{E}\big[\big|\int_{\tau}^{T}Z(s)dW(s)\big|^{2}\big|\mathcal{F}_{\tau}\big]=
\sup_{\tau}\mathbb{E}\big[\int_{\tau}^{T}Z^{2}(s)ds\big|\mathcal{F}_{\tau}\big]<\infty,\quad
\mathbb{P}-a.s.,
\end{eqnarray*}
where the supremum is taken over all stopping times $\tau$ with
respect to $\mathbb{F}$ and bounded by $T$, and $\mathcal{F}_\tau$
denotes the $\sigma$-algebra of the $\tau$-past.

\begin{lem}
Assume that $\xi$ is $\mathbb{P}$-a.s. bounded and \textbf{(A5)}
holds. The backward stochastic differential equation \eqref{bsdelin}
has a unique solution $(Y,Z)\in\mathbb{S}^{2}(\mathbb{R})\times
\mathbb{H}^{2}(\mathbb{R})$, and the integral process
$(\int_{0}^{t}Z(s)dW(s))_{0\leq t\leq T}$ is a BMO martingale.
\end{lem}
\Proof Take two stopping times $\tau_{1}, \tau_{2}\le T$ such that
$0\le \tau_{2}-\tau_{1}\leq \delta$ holds $\mathbb{P}$-a.s., with a
sufficiently small constant $\delta$ to be specified later in the
proof. The martingale representation \eqref{repr} gives
\begin{eqnarray*}
\mathbb{E}^{\mathbb{Q}}\big[\xi\big|\mathcal{F}_{\tau_{2}}\big]=
\mathbb{E}^{\mathbb{Q}}\big[\xi\big|\mathcal{F}_{\tau_{1}}\big]+\int_{\tau_{1}}^{\tau_{2}}Z(s)dW(s)-
\int_{\tau_{1}}^{\tau_{2}}Z(s)\alpha((s-T,0])g(s)ds.
\end{eqnarray*}
Now take a stopping time $\theta$ such that $\theta\leq \tau_{1}$
holds $\mathbb{P}$-a.s. As $Z\in\mathbb{H}^{2}(\mathbb{R})$, we have
\begin{eqnarray*}
\lefteqn{\mathbb{E}\big[\int_{\tau_{1}}^{\tau_{2}}|Z(s)|^{2}ds\big|\mathcal{F}_{\theta}\big]=
\mathbb{E}\big[\big|\int_{\tau_{1}}^{\tau_{2}}Z(s)dW(s)\big|^{2}\big|\mathcal{F}_{\theta}\big]}\nonumber\\
&=&\mathbb{E}\big[\big|\mathbb{E}^{\mathbb{Q}}\big[\xi\big|\mathcal{F}_{\tau_{2}}\big]
-\mathbb{E}^{\mathbb{Q}}\big[\xi\big|\mathcal{F}_{\tau_{1}}\big]+
\int_{\tau_{1}}^{\tau_{2}}Z(s)\alpha((s-T,0])g(s)ds\big|^{2}\big]\nonumber\\
&\leq&4\big(C+\delta
G\mathbb{E}\big[\int_{\tau_{1}}^{\tau_{2}}|Z(s)|^{2}ds\big|\mathcal{F}_{\theta}\big]\big),\quad
\mathbb{P}-a.s.,
\end{eqnarray*}
where $C$ denotes the uniform bound of $\xi$ and $G$ the one of $g$.
For a sufficiently small $\delta<\frac{1}{4G}$ we have
\begin{eqnarray*}
\mathbb{E}\big[\int_{\tau_{1}}^{\tau_{2}}|Z(s)|^{2}ds\big|\mathcal{F}_{\theta}\big]\leq
M, \quad \mathbb{P}-a.s.
\end{eqnarray*}
holds with a finite constant $M$ (independent of the stopping
times).\\
\noindent Now let $\delta<\frac{1}{4G}$, and take an arbitrary
stopping time $\tau$. Define $\tau_k = (\tau + k\delta)\wedge T,
k\ge 0.$ Then $(\tau_k)_{k\ge 0}$ is a sequence of stopping times
with respect to $\mathbb{F}$ such that $\tau_k - \tau_{k-1} \le
\delta,$ and $\tau_k-\tau_{k-1} = 0$ if $ k\ge N =
[\frac{T}{\delta}]+1.$ We can deduce from the inequality proved
before
\begin{eqnarray*}
\mathbb{E}\big[\int_{\tau}^{T}|Z(s)|^{2}ds\big|\mathcal{F}_{\tau}\big]=
\sum_{i=1}^{N}\mathbb{E}\big[\int_{\tau_{i-1}}^{\tau_i}|Z(s)|^{2}ds\big|\mathcal{F}_{\tau}\big]\le
N\, M,\quad \mathbb{P}-a.s.
\end{eqnarray*}
This proves the \emph{BMO} property.\cbdu

For generators without time delay, it is well known that the
solution component $Y$ inherits uniform boundedness from the
terminal condition $\xi$, see Proposition 2.1 in \cite{Kob}. We
shall now exhibit an example showing that this is not the case if
the generator has a linear delay dependence, as in Corollary
\ref{propint}.

{\bf Example 4}\\
Consider first the local martingale $M(t) = \int_0^t
\frac{2}{(1-s)^3} d W^\mathbb{Q}(s),\,\, t\in[0,1)$, under the
measure $\mathbb{Q}$ defined in Corollary \ref{propint}. Let
$$\tau = \inf \{t\ge 0: |M(t)| \ge 1\} \wedge 1.$$
We start by showing that $\tau$ can take values arbitrarily close to
1 with positive probability for $\mathbb{Q}$, hence also for the
equivalent $\mathbb{P}.$ This claim follows from classical results.
In fact, by time change with the quadratic variation
$$\langle M\rangle_t = \int_0^t \frac{4}{(1-s)^6} ds =
\frac{4 q(t)}{5(1-t)^5},\quad t\in[0,1),$$ with $q(t) = 1 -(1-t)^5,
\,\,t\in[0,1),$ $(M(t))_{0\leq t<1}$ has the same law as
$(B(\frac{4q(t)}{5(1-t)^5}))_{0\leq t<1}$ with a
$\mathbb{Q}$-Brownian motion $B$. Defining
$$\sigma = \inf\{ t\ge 0: |B(t)|\ge 1\},$$
we obtain that under $\mathbb{Q}$, $\sigma$ has the same law as
$\frac{4 q(\tau)}{5(1-\tau)^5}.$ Since $\sigma$ can take values
arbitrarily close to $\infty$ with positive probability, $\tau$ can
take values arbitrarily close to $1$ with positive probability.\\
\indent Consider the linear BSDE
\begin{eqnarray*}
Y(t)=\xi+\int_{t}^{1}\int_{0}^{s}Z(u)du
ds-\int_{t}^{1}Z(s)dW(s),\quad 0\leq t\leq 1,
\end{eqnarray*}
and define $\xi= M_\tau.$ Then $\xi$ is bounded, and we have the
stochastic integral representation
$$\xi = \int_0^1 \frac{2}{(1-s)^3} 1_{[0,\tau)}(s) d W^\mathbb{Q}(s).$$
Therefore, Corollary \ref{propint} yields the solution
\begin{eqnarray*}
Z(t) &=& \frac{2}{(1-t)^3} 1_{[0,\tau)}(t),\quad 0\leq t\leq 1,\\
Y(t)&=&\mathbb{E}^{\mathbb{Q}}[M_{\tau}|\mathcal{F}_{t}]+(1-t)\int_{0}^{t}Z(s)ds\\
&=&\mathbb{E}^{\mathbb{Q}}[M_{\tau}|\mathcal{F}_{t}]+(1-t)\frac{(\tau\wedge
t)(2-\tau\wedge t)}{1-(\tau\wedge t)^{2}} ,\quad 0\leq t\leq 1.
\end{eqnarray*}
Take an arbitrary constant $C>0$. We can find $u\in[0,1)$ such that
$\frac{u(2-u)}{1-u}>C-1$. As $\tau$ takes values arbitrarily close
to 1 with positive probability, $\mathbb{P}(\tau\geq u)>0$ and with
positive probability
$$(1-u)\frac{(\tau\wedge
u)(2-\tau\wedge u)}{1-(\tau\wedge u)^{2}}=\frac{u(2-u)}{1-u}>C-1,$$
$$Y(u)>C.$$
We can conclude that for an arbitrarily large $C$ there exists
$u\in[0,1)$ such that the process $Y$ at time $u$ crosses $C$ with
positive probability, $\mathbb{P}(Y(u)>C)>0$. This establishes the
lack of uniform boundedness for $Y$.

\section{The comparison principle and measure solutions}
\indent The concepts of comparison principle and measure solutions
play an important role in the theory of BSDE without time delays. In
this section we first show by an example that they cannot be
extended to time delayed BSDE. We shall see that the failure of the
properties can be traced back to a sign change of the control
process $Z$, and consequently show that they continue to hold on
stochastic intervals on which $Z$ stays away from 0. For a statement
of the comparison principle, we refer the reader to Theorem 2.2 in
\cite{K},
and for the concept of a measure solution to the paper \cite{I2}.\\
\indent {\bf Example 5}\\
We first give an example exhibiting a failure of the comparison
principle. Consider the linear backward stochastic differential
equation with time delayed generator
\begin{eqnarray}\label{exam11}
Y(t)=\xi+\int_{t}^{T}\int_{0}^{s}Z(u)du
ds-\int_{t}^{T}Z(s)dW(s),\quad 0\leq t\leq T.
\end{eqnarray}
Take $\tilde{\xi}=0$. The corresponding solution of \eqref{exam11}
is $\tilde{Y}(t)=\tilde{Z}(t)=0,0\leq t\leq T$. To compare, take
$\xi=\big(W(T)-\frac{T^{2}}{2}\big)^{2}$. By applying Corollary
\ref{propint} we can construct the corresponding solution $(Y,Z)$ of
\eqref{exam11}. The martingale representation of $\xi$ under the
measure $\mathbb{Q}$ with the $\mathbb{Q}$-Brownian motion
$W^{\mathbb{Q}}(t) = W(t) + \frac{(T-t)^2}{2} - \frac{T^2}{2}, 0\le
t\le T,$ yields
\begin{eqnarray*}
\mathbb{E}^{\mathbb{Q}}[\big(W(T)-\frac{T^{2}}{2}\big)^{2}|\mathcal{F}_{t}]
&=&\mathbb{E}^{\mathbb{Q}}[W^{\mathbb{Q}}(T)^{2}|\mathcal{F}_{t}]\\
&=& W^{\mathbb{Q}}(t)^{2}+(T-t)\nonumber\\
&=&T+\int_{0}^{t}2W^{\mathbb{Q}}(s)\,dW^{\mathbb{Q}}(s),
\end{eqnarray*}
and we can derive the solution
\begin{eqnarray*}
Y(t)=W^{\mathbb{Q}}(t)^{2}+(T-t)+(T-t)\int_{0}^{t}2W^{\mathbb{Q}}(t)ds.
\end{eqnarray*}
Clearly, $\tilde{\xi}\leq \xi$ holds $\mathbb{P}$-a.s. It is
straightforward to note that $\tilde{Y}(0)\leq Y(0)$ holds
$\mathbb{P}$-a.s. However, we claim that for any $t\in(0,T)$ we have
$\mathbb{Q}\big(\tilde{Y}(t)>Y(t)\big)>0$, which, by equivalence of
the measures, contradicts the comparison principle under
$\mathbb{P}$.\\
\indent To prove that $\mathbb{Q}\big(\tilde{Y}(t)>Y(t)\big)>0$, it
is sufficient to show that the conditional law of
$\int_{0}^{t}W^{\mathbb{Q}}(s)ds$ given $W^{\mathbb{Q}}(t)=x$ is
unbounded for any $t\in(0,T)$ and any $x\in\mathbb{R}$. This can be
verified under $\mathbb{P}$ for $W$ instead of $W^{\mathbb{Q}}$ as
well. First recall, see Chapter 5.6.B in \cite{Kar}, that the
conditional law of $W(s)$ given $W(t)=x$ is nondegenerate Gaussian,
for $0<s<t$, for any $(t,x)\in(0,T)\times\mathbb{R}$, and that the
process $[0,t]\ni s\mapsto W(s)$ is a Brownian bridge from $0$ to
$x$ conditional on $\{W_t = x\}$. We have the convergence
\begin{eqnarray*}
\int_{0}^{t}W(s)ds=\lim_{N\rightarrow\infty}\frac{\sum_{i=1}^{N}W(t_{i}^N)}{N},\quad
\mathbb{P}-a.s.\quad \mbox{and in}\quad L^2,
\end{eqnarray*}
for a sequence of equidistant partitions
$0=t_{0}<t_{1}^N<...<t_{N}^N=t$ of $[0,t]$. As the $L^2$-limit of a
Gaussian sequence $\int_0^t W(s) ds$ is Gaussian. It is
straightforward to show that the variance of $\int_{0}^{t}W(s)ds$
conditional on $\{W_t = x\}$ is strictly positive. We conclude that
the conditional law of $\int_{0}^{t}W(s)ds$ given
$\{W(t)=x\}$ is unbounded for any $t\in(0,T)$ and any $x\in\mathbb{R}$.\\

\indent A failure of the comparison principle indicates that it may
also not always be possible to represent a solution of a BSDE with a
time delayed generator as a measure solution. Recall, that for a
BSDE without time delay and with a Lipschitz continuous generator
independent of $Y$, a unique square integrable solution can always
be represented as a conditional expectation of the terminal value
under an appropriate probability measure (a \emph{measure
solution}), see Theorem 1.1 in \cite{I2}. The following example
shows that this property
may fail for the solution of a time delayed BSDE.  \\

\indent {\bf Example 6}\\ Consider again the linear backward
stochastic differential equation
\begin{eqnarray}\label{ms}
Y(t)=\xi+\int_{t}^{T}\int_{0}^{s}Z(u)du
ds-\int_{t}^{T}Z(s)dW(s),\quad 0\leq t\leq T,
\end{eqnarray}
with $\xi=\int_{0}^{T}\cos tdW(t)$. An easy calculation shows that
there exists a unique square integrable solution given by
\begin{eqnarray*}
Z(t)&=&\cos t,\quad 0\leq t\leq T,\nonumber\\
Y(t)&=&\int_{0}^{t}\cos sdW(s)+\cos t-\cos T,\quad 0\leq t\leq T.
\end{eqnarray*}
To describe a possible measure solution, for $T<\frac{\pi}{2}$
rewrite the equation \eqref{ms} as
\begin{eqnarray*}
Y(t)=\xi+\int_{t}^{T}\cos s \big(dW(s)-\tan sds\big),\quad 0\leq
t\leq T.
\end{eqnarray*}
One can define the equivalent probability measure
\begin{eqnarray*}
\frac{d\mathbb{Q}}{d\mathbb{P}}\big|_{\mathcal{F}_{T}}=\exp\big(\int_{0}^{T}\tan
sdW(s)-\frac{1}{2}\int_{0}^{T}\tan^{2}sds\big),
\end{eqnarray*}
and the unique measure solution under this measure
\begin{eqnarray*}
Y(t)=\mathbb{E}^{\mathbb{Q}}\big[\xi|\mathcal{F}_{t}\big],\quad
0\leq t\leq T.
\end{eqnarray*}
Consider now the case $T=\frac{\pi}{2}$. If there were a measure
solution on $[0,\frac{\pi}{2}]$ under some probability measure
$\mathbb{Q}$ then
\begin{eqnarray*}
\frac{d\mathbb{Q}}{d\mathbb{P}}\big|_{\mathcal{F}_{t}}=\exp\big(\int_{0}^{t}\tan
sdW(s)-\frac{1}{2}\int_{0}^{t}\tan^{2}sds\big),
\end{eqnarray*}
for any $0\leq t<\frac{\pi}{2}$, and the following limiting relation
would hold
\begin{eqnarray*}
\lim_{t\rightarrow \frac{\pi}{2}}\exp\big(\int_{0}^{t}\tan s
dW(s)-\frac{1}{2}\int_{0}^{t}\tan^{2}sds\big)=\frac{d\mathbb{Q}}{d\mathbb{P}}\big|_{\mathcal{F}_{\frac{\pi}{2}}}.
\end{eqnarray*}
We show that this limit is not a probability density.\\
\noindent Define a sequence of points
$0=t_{0}<t_{1}<...<t_{n}<\cdots<\frac{\pi}{2}$, for
$n\in\mathbb{N}$, such that
\begin{eqnarray*}
\int_{t_{i-1}}^{t_{i}}\tan^{2} sds=1,\quad \forall i\in\mathbb{N},
\end{eqnarray*}
and $\lim_{n\rightarrow \infty} t_n=\frac{\pi}{2}$. The sequence of
random variables $(X_{i})_{i\in\mathbb{N}}$ defined by
$X_{i}=\int_{t_{i-1}}^{t_{i}}\tan sdW(s)$ is i.i.d. with standard
Gaussian laws. The strong law of large numbers implies
\begin{eqnarray*}
\lefteqn{\lim_{t\to \frac{\pi}{2}}\exp\big(\int_{0}^{t}\tan
s dW(s)-\frac{1}{2}\int_{0}^{t}\tan^{2}sds\big)}\\
&=&\lim_{n\to \infty}\exp\big(\int_{0}^{t_n}\tan
s dW(s)-\frac{1}{2}\int_{0}^{t_n}\tan^{2}sds\big)\nonumber\\
&=&\lim_{n\to
\infty}\exp\big(X_{1}+\cdots + X_{n}-\frac{1}{2}n\big)\nonumber\\
&=& \lim_{n\to
\infty}\exp\big(n(\frac{X_{1}+\cdots + X_{n}}{n}-\frac{1}{2})\big)\nonumber\\
&=&0,\quad \mathbb{P}-a.s.
\end{eqnarray*}
This shows that an equivalent probability measure cannot be defined
on $[0,\frac{\pi}{2}]$. In summary we have established a BSDE with a
time delayed generator, for which there is a unique square
integrable solution, whereas a measure solution fails to exist.\\

\indent One observation we can draw from the preceding two examples
is that the comparison principle may not hold and measure solutions
may fail to exist, while the (continuous) control process $Z$ can
cross $0$. In the following two Theorems we shall exclude the
approach of the difference of two control processes resp. one
control process to 0 or $\infty$ by stopping them before passages of
small or large thresholds happen. We shall prove that on the
corresponding stochastic intervals, the comparison principle holds,
and a measure solution exists.
\begin{thm}
Consider the backward stochastic differential equations \eqref{bsde}
with generators $f, \tilde{f}$ and corresponding terminal values
$\xi, \tilde{\xi}$ satisfying the assumptions
\textbf{(A1)}-\textbf{(A4)}. Let $(Y,Z)$ and $(\tilde{Y},\tilde{Z})$
denote the associated unique solutions in
$\mathbb{S}^{2}(\mathbb{R})\times \mathbb{H}^{2}(\mathbb{R})$. For
$n\in\mathbb{N}$ define the stopping time $\tau_{n}=\inf\{t\geq 0:
|Y(t)-\tilde{Y}(t)|\vee|Z(t)-\tilde{Z}(t)|\leq \frac{1}{n}$ or
$|Y(t)-\tilde{Y}(t)|\vee |Z(t)-\tilde{Z}(t)|\geq n\}\wedge T$.
Suppose that
\begin{itemize}
\item $Y(\tau_{n})\geq \tilde{Y}(\tau_{n}), \mathbb{P}$-a.s.,
\item $\delta f(t,y_{t},z_{t}):=f(t,y_{t},z_{t})- \tilde{f}(t,y_{t},z_{t})\geq 0$,\quad $t\in[0,T], (y_t, z_t)\in L_{-T}^\infty(\mathbb{R})\times L_{-T}^2(\mathbb{R})$.
\end{itemize}
Then, $Y(t)\geq \tilde{Y}(t)$ holds $\mathbb{P}$-a.s on
$[0,\tau_{n}]$.
\end{thm}
\Proof We follow the idea from the proof of Theorem 2.2 in \cite{K}.
For $t\in[0,T]$ let
$$\delta Y(t)=Y(t)-\tilde{Y}(t),\quad \delta Z(t)= Z(t)-\tilde{Z}(t),$$
$$\delta f(s,\tilde{Y}_{s},\tilde{Z}_{s}) = f(s, \tilde{Y}_{s},\tilde{Z}_{s})-\tilde{f}(s,\tilde{Y}_{s},\tilde{Z}_{s}),$$
and
$$\triangle_{y}f(t)=\frac{f(t,Y_{t},Z_{t})-
f(t,\tilde{Y}_{t},Z_{t})}{Y(t)-\tilde{Y}(t)},\quad
\triangle_{z}f(t)=\frac{f(t,\tilde{Y}_{t},Z_{t})-f(t,\tilde{Y}_{t},\tilde{Z}_{t})}{Z(t)-\tilde{Z}(t)}.$$
We can derive
\begin{eqnarray*}
\delta Y(t)&=&\delta Y(\tau_{n})+\int_{t}^{\tau_{n}}\big(\triangle_{y}f(s)\delta Y(s)+\triangle_{z}f(s)\delta Z(s)+\delta f(s,\tilde{Y}_{s},\tilde{Z}_{s})\big)ds\nonumber\\
&&-\int_{t}^{T}\delta Z(s)dW(s),\quad 0\leq t\leq \tau_{n}.
\end{eqnarray*}
By rewriting these expressions and changing the measure we obtain
\begin{eqnarray}\label{delta}
\delta Y(t)&=&\delta Y(\tau_{n})e^{\int_{t}^{\tau_{n}}\triangle_{y}f(s)ds}+\int_{t}^{\tau_{n}}\delta f(s,\tilde{Y}_{s},\tilde{Z}_{s})e^{\int_{t}^{s}\triangle_{y}f(u)du}ds\nonumber\\
&&-\int_{t}^{\tau_{n}}\delta
Z(s)e^{\int_{t}^{s}\triangle_{y}f(u)du}dW^{\mathbb{Q}}(s),\quad
0\leq t\leq \tau_{n},
\end{eqnarray}
with the equivalent probability measure $\mathbb{Q}$ defined as
\begin{eqnarray*}
\frac{d\mathbb{Q}}{d\mathbb{P}}\big|_{\mathcal{F}_{\tau_{n}}}=\exp\big(\int_{0}^{\tau_{n}}
\triangle_{z}f(s)dW(s)-\frac{1}{2}\int_{0}^{\tau_{n}}\big(\triangle_{z}f(s)\big)^{2}ds\big).
\end{eqnarray*}
Since $\delta Z\in\mathbb{H}^{2}(\mathbb{R})$ under the measure
$\mathbb{P}$, the density
$\frac{d\mathbb{Q}}{d\mathbb{P}}\big|_{\mathcal{F}_{\tau_{n}}}$ is
square integrable under the measure $\mathbb{P}$, and
$t\mapsto\triangle_{y}f(t)$ is a.s. uniformly bounded up to time
$\tau_{n}$, we can use Cauchy-Schwarz' inequality to obtain
\begin{eqnarray*}
\mathbb{E}^{\mathbb{Q}}\big[\sqrt{\int_{0}^{\tau_{n}}\big|\delta
Z(s)e^{\int_{0}^{s}\triangle_{y}f^{1}(u)du}\big|^{2}ds}\big]<\infty,
\end{eqnarray*}
and $\big(\int_{0}^{t}\delta
Z(s)e^{\int_{0}^{s}\triangle_{y}f^{1}(u)du}dW^{\mathbb{Q}}(s)\big)_{0\leq
t\leq \tau_{n}}$ is a $\mathbb{Q}$-martingale with vanishing
expectation, see Theorem 3.28 in \cite{Kar}. Taking the conditional
expectation with respect to $\mathcal{F}_t$ on both sides of
equation \eqref{delta} under the measure $\mathbb{Q}$, we get the
desired result. \cbdu

\begin{thm}\label{thms}
Consider the backward stochastic differential equation \eqref{bsde}
with the generator $f(t,y_{t},z_{t})=f(t,z_{t}), t\ge 0, (y_t,
z_t)\in L_{-T}^\infty(\mathbb{R})\times L_{-T}^2(\mathbb{R})$, and
the corresponding terminal value $\xi$ satisfying the assumptions
\textbf{(A1)}-\textbf{(A4)}. Let $(Y,Z)$ denote the associated
unique solution in $\mathbb{S}^{2}(\mathbb{R})\times
\mathbb{H}^{2}(\mathbb{R})$. For $n\in\mathbb{N}$ define the
stopping time $\tau_{n}=\inf\{t\geq 0; |Z(t)|\leq
\frac{1}{n}\,\,\mbox{or}\,\, |Z(t)|\geq n\}\wedge T$. Then, there
exists a unique equivalent probability measure $\mathbb{Q}$
restricted to $[0,\tau_{n}]$ such that
\begin{eqnarray*}
Y(t)=\mathbb{E}^{\mathbb{Q}}\big[Y(\tau_{n})|\mathcal{F}_{t}\big],\quad
\end{eqnarray*}
holds for all $t\in[0,\tau_{n}]$, $\mathbb{P}$-a.s.
\end{thm}
\Proof The proof requires a change of measure argument, just as in
the preceding proof. Details are omitted. \cbdu

\indent We remark that if $Y$ is $\mathbb{P}$-a.s. bounded, one can
define $\tau^{m}_{n}=\inf\{t\geq 0: |Z(t)|\leq
\frac{1}{n}\,\,\mbox{\emph{or}}\,\, |Z(t)|\geq m\}\wedge T$ for
$n,m\in\mathbb{N}$ and show that the corresponding family of
measures
$\big(\frac{d\mathbb{Q}}{d\mathbb{P}}\big|_{\mathcal{F}_{\tau^{m}_{n}}}\big)_{m\in\mathbb{N}}$
is uniformly integrable. Compare the proof of Theorem 1.1 \cite{I2}.
In this case a unique measure solution can be defined on
$[0,\tau^{\infty}_{n}]$
with $\tau^{\infty}_{n}=\inf\{t\geq 0; |Z(t)|\leq \frac{1}{n}\}\wedge T$.\\
\
\\
\
\\
\noindent \textbf{Acknowledgements:} This paper was written while
the first author was staying at Humboldt University Berlin.
{\L}ukasz Delong acknowledges the financial support from AMaMeF
programme.


\begin{thebibliography}{99}

\bibitem{I2}
Ankirchner, S., Imkeller, P., Popier, A. (2008) On measure solutions
of backward stochastic differential equations. Preprint.
\bibitem{Bech}
Becherer, D. (2006) Bounded solutions to backward SDE's with jumps
for utility optimization and indifference pricing. Ann. Appl.
Probab. \textbf{16}, 2027-2054.
\bibitem{BS}
Blanchet-Scalliet, C., Eyraud-Loisel, A., Royer-Carenzi, M. (2008).
Hedging of defaultable claims using BSDE with uncertian time
horizon. Preprint.
\bibitem{BI}
Buckdahn, R., Imkeller, P. (2009) Backward stochastic differential
equations with time delayed generator. Preprint.
\bibitem{DI}
Delong, {\L}., Imkeller, P. (2009) On Malliavin's differentiability
of time delayed BSDEs driven by Brownian motions and Poisson random
measures. Preprint.
\bibitem{K}
El Karoui, N., Peng, S., Quenez, M.C., (1997) Backward stochastic
differential equations in finance. Math. Finance \textbf{7}, 1-71.
\bibitem{FT}
Fuhrman, M., Masiero, F., Tessitore, G. (2008) Stochastic equations
with delay: optimal control via BSDEs and regular solutions of
Hamilton-Jacobi-Bellman equations. Preprint.
\bibitem{I1}
Imkeller, P. (2008) \emph{Malliavin's Calculus and Applications in
Stochastic Control and Finance.} IM PAN Lectures Notes, Warsaw.
\bibitem{Kar}
Karatzas, I., Shreve, S. (1988) \emph{Brownian Motion and Stochastic
Calculus}. Springer-Verlag.
\bibitem{Kob}
Kobylanski, M. (2000) Backward stochastic differential equations and
partial differential equations with quadratic growth. Ann. Probab.
\textbf{28}, 558-602.
\bibitem{Moh}
Mohammed, S.E.A. (1984) \emph{Stochastic functional Differential
Equations}. Pitman.
\bibitem{Mor}
Morlais, M. A. (2009) Quadratic BSDEs driven by a continuous
martingale and applications to the utility maximization problem.
Finance Stoch. \textbf{13}, 121-150.
\bibitem{PP}
Pardoux, E. Peng, S. (1990) Adapted solution of a backward
stochastic differential equation. Systems Control Lett. \textbf{14},
55-61.
\bibitem{Pr}
Protter, P. (1992) \emph{Stochastic Integration and Differential
Equations}. Springer-Verlag.
\bibitem{Xu}
Xu, D., Yang, Z., Huang,Y. (2008) Existence-uniqueness and
continuation theorems for stochastic functional differential
equations. J. Differential Equations \textbf{254}, 1681-1703.
\end{thebibliography}
\end{document}